\documentclass[11pt]{article}
\pdfoutput=1
\usepackage{amscd,amsmath, amssymb, fancyhdr, epsfig, color, tocloft, mathtools}
\usepackage{dutchcal}

\usepackage[backref=page]{hyperref}   
\renewcommand*{\backrefalt}[4]{%
	\ifcase #1 (Not cited.)%
	\or        (Cited on page~#2.)%
	\else      (Cited on pages~#2.)%
	\fi}

\hypersetup{
	colorlinks   = true,
	citecolor    = magenta,
	linkcolor    = blue,
	urlcolor     = magenta	
}

\newcommand{\version}{version 2, October 13, 2023}

\numberwithin{equation}{section}

\newcommand{\RR}{\mathbb{R}}

\newcommand{\GL}{{\sf{GL}}}

\renewcommand{\O}{{\sf{O}}}
\newcommand{\CO}{{{CO}}}
\renewcommand{\H}{{\sf{H}}}

\def\eqref#1{(\ref{#1})}

\newcommand{\goth}{\mathfrak}

\newcommand{\ra}{{\:\longrightarrow\:}}
\newcommand{\N}{{\mathbb N}}
\newcommand{\Z}{{\mathbb Z}}
\newcommand{\C}{{\mathbb C}}
\newcommand{\R}{{\mathbb R}}

\renewcommand{\H}{{\mathbb H}}

\def\1{\sqrt{-1}\,}

\newcommand{\restrict}[1]{{\left|_{{\phantom{|}\!\!}_{#1}}\right.}}

\newcommand{\cntrct}                
{\hspace{2pt}\raisebox{1pt}{\text{$\lrcorner$}}\hspace{2pt}}

\newcommand{\arrow}{{\:\longrightarrow\:}}

\renewcommand{\bar}{\overline}
\renewcommand{\phi}{\varphi}
\renewcommand{\epsilon}{\varepsilon}
\renewcommand{\geq}{\geqslant}

\newcommand{\Hom}{\operatorname{Hom}}

\renewcommand{\dim}{\operatorname{\sf dim}}

\renewcommand{\Im}{\operatorname{Im}}

\newcounter{Mycounter}[section]
\newcounter{lemma}[section]
\renewcommand{\thelemma}{{Lemma \thesection.\arabic{lemma}}}
\newcommand{\lemma}{%
	\setcounter{lemma}{\value{Mycounter}}
	\refstepcounter{lemma}
	\stepcounter{Mycounter}
	{\noindent \bf \thelemma:\ }}

\newcounter{claim}[section]
\newcounter{sublemma}[section]
\newcounter{corollary}[section]
\renewcommand{\thecorollary}{{Corollary \thesection.\arabic{corollary}}}
\newcommand{\corollary}{%
	\setcounter{corollary}{\value{Mycounter}}
	\refstepcounter{corollary}
	\stepcounter{Mycounter}
	{\noindent \bf \thecorollary:\ }}

\newcounter{theorem}[section]
\renewcommand{\thetheorem}{{Theorem \thesection.\arabic{theorem}}}
\newcommand{\theorem}{%
	\setcounter{theorem}{\value{Mycounter}}
	\refstepcounter{theorem}
	\stepcounter{Mycounter}
	{\noindent \bf \thetheorem:\ }}

\newcounter{conjecture}[section]
\newcounter{proposition}[section]
\renewcommand{\theproposition} {{Proposition \thesection.\arabic{proposition}}}
\newcommand{\proposition}{%
	\setcounter{proposition}{\value{Mycounter}}
	\refstepcounter{proposition}
	\stepcounter{Mycounter}
	{\noindent \bf \theproposition:\ }}

\newcounter{definition}[section]
\renewcommand{\thedefinition} {{Definition~\thesection.\arabic{definition}}}
\newcommand{\definition}{%
	\setcounter{definition}{\value{Mycounter}}
	\refstepcounter{definition}
	\stepcounter{Mycounter}
	{\noindent \bf \thedefinition:\ }}

\newcounter{example}[section]
\renewcommand{\theexample}{{Example \thesection.\arabic{example}}}
\newcommand{\example}{%
	\setcounter{example}{\value{Mycounter}}
	\refstepcounter{example}
	\stepcounter{Mycounter}
	{\noindent \bf \theexample:\ }}

\newcounter{remark}[section]
\renewcommand{\theremark}{{Remark \thesection.\arabic{remark}}}
\newcommand{\remark}{%
	\setcounter{remark}{\value{Mycounter}}
	\refstepcounter{remark}
	\stepcounter{Mycounter}
	{\noindent \bf \theremark:\ }}

\newcounter{problem}[section]
\newcounter{question}[section]
\makeatletter

\def\blacksquare{\hbox{\vrule width 5pt height 5pt depth 0pt}}
\def\endproof{\blacksquare}

\newcommand{\proof}{{\bf Proof: \ }}
\newcommand{\pstep}{{\bf Proof. Step 1: \ }}

\begin{document}
	
	\begin{center}
		{\Large\bf  Do products of compact complex manifolds\\[.1in] admit LCK metrics?}\\[5mm]
		{\large
			Liviu Ornea\footnote{Partially supported by Romanian Ministry of Education and Research, Program PN-III, Project number PN-III-P4-ID-PCE-2020-0025, Contract  30/04.02.2021.}, Misha Verbitsky\footnote{Partially supported by
				by the HSE University Basic Research Program, FAPERJ E-26/202.912/2018 
				and CNPq - Process 310952/2021-2.}, 
			Victor Vuletescu\footnote{Partially supported by Romanian Ministry of Education and Research, Program PN-III, Project number PN-III-P4-ID-PCE-2020-0025, Contract  30/04.02.2021.
			\\[1mm]
				\noindent{\bf Keywords:} LCK manifold, product manifold, Hopf surface, Inoue surface, rational curve. 
				
				\noindent {\bf 2020 Mathematics Subject Classification:} {32J27, 53C55}
			}\\[4mm]
			
		}
		
	\end{center}

	{\small
		\hspace{0.15\linewidth}
		\begin{minipage}[t]{0.7\linewidth}
			{\bf Abstract} \\ 
An LCK (locally conformally K\"ahler) manifold is a
Hermitian manifold which admits a K\"ahler cover with deck
group acting by holomorphic homotheties with respect to
the K\"ahler metric. The product of two LCK manifolds does
not have a natural product LCK structure. It is
conjectured that a product of two compact complex
manifolds is never LCK. We classify all known examples of
compact LCK manifolds onto three not exclusive classes: LCK
with potential, a class of manifolds we call  of Inoue type, and those containing a
rational curve. In the present paper, we prove that a
product of an LCK manifold and an LCK manifold belonging
to one of these three classes does not admit an LCK
structure.
		\end{minipage}
	}

	\tableofcontents
	
	
\section{Introduction}

A {\bf locally conformally K\"ahler} (LCK) manifold is a Hermitian manifold $(M,I,g)$ which admits a K\"ahler cover such that the deck group acts by holomorphic homotheties with respect to the K\"ahler metric. Equivalently,  the fundamental form $\omega(x,y):=g(Ix,y)$, $x,y\in TM$, satisfies the integrability condition $d\omega=\theta\wedge\omega$ for some {\em closed} 1-form $\theta$, called {\bf the Lee form}. 

Typical examples are the Hopf manifolds, the Kodaira
surfaces, the Kato  manifolds. We refer to the book
\cite{_OV_book_} for details about LCK geometry. In
Section \ref{_Basics_and_examples_Section_} of this paper
we present the necessary background of LCK geometry
(Subsection \ref{_Definitions_and_basics_Subsection_})
and we describe all the known classes of examples of LCK manifolds
(Subsections \ref{_LCK_pot_Subsection_},
\ref{_LCK_Inoue_type_Subsection_}, \ref{_LCK
  with_rational_curves_Subsection_}).

One can easily see that if $(M,I,g)$ is LCK with Lee form $\theta$, then $(M,I,e^fg)$ is LCK with Lee form $\theta+df$. As such, LCK geometry is a part of conformal geometry. Since a product of conformal groups is not a conformal group ($\CO(n)\times \CO(m)\not\subset\CO(n+m)$), the product of two LCK manifolds does not have a natural product LCK structure. 

Note that assuming from the beginning that the factors of the product are LCK is not restrictive. Indeed, if the product admits an LCK structure, then both factors have induced LCK structures, because a complex submanifold of an LCK manifold is again LCK.

The question then arises whether the product of two LCK
manifolds does admit any LCK structure compatible with the
product complex structure. Even more generally (see
\cite[Question 47.6]{_OV_book_}), does the product of two
complex manifolds admit any LCK structure, not necessarily
compatible with the product complex structure? Only
partial answers to this question were given up to now. The following products are known to not bear strict LCK metrics: the product of two compact Vaisman manifolds (\cite[Corollary 3.3]{tsu}), the product of a compact K\"ahler manifold of dimension at least 2 and a compact strict LCK manifold (\cite[Corollary 2]{opv}), the product of a compact complex smooth curve with a compact manifold carrying no LCK metric with potential (\cite[Proposition 7]{nico2}). 

We observe that all known examples of compact LCK manifolds fall in one of the following three classes: (1) LCK manifolds with potential, (2) LCK manifolds that we call {\bf of Inoue type}, a class containing all LCK Inoue surfaces and the LCK Oeljeklaus-Toma manifolds, and (3) LCK manifolds containing rational curves. We define these classes and provide appropriate examples in Subection \ref{_LCK_pot_Subsection_},  Subection \ref{_LCK_Inoue_type_Subsection_}, and Subection \ref{_LCK with_rational_curves_Subsection_} respectively. Note that these three classes are not mutually exclusive, see \ref{_blow_up_of_Inoue_type_Remark_}.


The goal of this paper is to show that if $X$ is a compact LCK manifold in any of these classes above and $Y$ is an arbitrary compact complex manifold, then the product $M:=X\times Y$ does not carry any LCK metric (see \ref{_nu_pot_}, \ref{_nu_inoue_} and \ref{_nu_spl_}). 

Note that if we admit the Global Spherical Shell
conjecture (Subsection \ref{_LCK with_rational_curves_Subsection_}), 
then the 3 classes above cover all LCK
surfaces. It follows that, if the GSS conjecture is true,
then  a product of a compact LCK surface and any
compact complex manifold is not of LCK type.

\section{Preliminaries}\label{_Basics_and_examples_Section_}

In this section, we gather the very basic definitions and results in locally conformally K\"ahler geometry. We refer to the recent monograph \cite{_OV_book_}.

\subsection{Definitions and basic results}\label{_Definitions_and_basics_Subsection_}

Let $(M,I,g,\omega)$ be a Hermitian manifold, $\dim_\C M\geq 2$. Here $\omega(\cdot, \cdot)=g(I\cdot, \cdot)$.

\hfill

\definition\label{_lck_}
The Hermitian manifold $(M,I,g,\omega)$ is {\bf locally conformally K\"ahler} (LCK) if 
there exists a closed 1-form $\theta$ such that $d\omega=\theta\wedge\omega$.
The 1-form $\theta$ is called the {\bf Lee form} and the 
$g$-dual vector field $\theta^\sharp$ is called the {\bf Lee field}.

If the Lee form is exact,  $(M,I,g,\omega)$ is called {\bf globally conformally K\"ahler (GCK)}. An LCK structure which is not GCK will be called {\bf strict}.

\hfill

\definition\label{_wlck_}
{More generally, if $(M, I)$ is a complex manifold, a {\bf weak lck structure (WLCK)} on $M$ is a $(1,1)-$form $\omega$ which obeys a relation $d\omega=\theta\wedge\omega$ for some closed 1-form $\theta$ and such that $h(\cdot, \cdot):=\omega(I\cdot, \cdot)$ is  positive definite outside a proper analytic subspace ${\mathcal B}(\omega)$ (called also ``the bad locus" of $\omega$).  A {\bf globally WLCK (GWLCK)} structure is one such that $\theta$ is exact.}

\hfill

\theorem\label{_LCK_cover_equiv_theorem_}(\cite[Remark, p. 236]{va_gd}, \cite[Section 3.4.2]{_OV_book_}) \\
The Hermitian manifold $(M,I,g,\omega)$ is LCK if and only if it admits a K\"ahler cover $\Gamma\arrow(\tilde M,\tilde\omega)\arrow M$ such that the deck group $\Gamma$ acts by holomorphic homotheties. \endproof

\hfill

It follows that for any LCK manifold $M$, one can define a
{\bf homothety character} $\chi:\Gamma \arrow \R^{>0}$
associating the scale
factor $\frac{\gamma^*\tilde\omega}{\tilde \omega}$
to each deck transform $\gamma$. The homothety character
is a representation of $\pi_1(M)$ and corresponds to the class of the Lee form
	$[\theta]$  under the canonical isomorphisms $ H^1(M,\R)\simeq \Hom_\Z\left(
	H_1(M), \R\right)\simeq
	{\Hom_\Z}\left(\frac{\pi_1(M)}{[\pi_1(M),\pi_1(M)]},
	\R\right)$. 

\hfill

\example\label{_class_Hopf_Example_}
The {\bf classical Hopf manifold}
is a quotient of $\C^n \backslash 0$ by
the $\Z$-action generated by a linear map
$z \mapsto \lambda z$, where $\lambda\in \C$,
$|\lambda| > 1$. This $\Z$-action 
multiplies the standard flat K\"ahler form
$-\1 \sum_i dz_i \wedge d\bar z_i$
by $|\lambda|^2$, hence defines
an LCK form $\frac{-\1 \sum_i dz_i \wedge d\bar
  z_i}{|z|^2}$ on the quotient 
$\frac {\C^n \backslash 0}{\Z}$.

\hfill

\remark\label{_conf_change_Remark_} The LCK condition is
conformally invariant: if $(M,I,g,\omega)$ is LCK with Lee
form $\theta$ and $f:M\rightarrow \R$ is smooth, then
$(M,I,e^fg,e^f\omega)$ is LCK with Lee form $\theta+df$.

\hfill

\remark Let $\iota:N\rightarrow M$ be a submanifold of an
LCK manifold $(M,I,g,\omega, \theta)$. Then $N$ is LCK
with induced complex structure and metric, and with Lee
form  $\iota^*\theta$. If $\iota^*\theta$ is exact, the
submanifold is called {\bf induced globally conformally
  K\"ahler (IGCK)}. We stress that not all K\"ahler
submanifolds of an LCK manifold are IGCK.

\hfill

\example
Let $M$ be a classical Hopf manifold, 
\ref{_class_Hopf_Example_}, and 
$E= \frac{\C\backslash 0}{\Z}\subset M$
an elliptic curve obtained from a complex
line in $\C^n\backslash 0$. Clearly, 
$E$ is K\"ahler, but the Lee form 
$\theta= - d \log|z|$ is clearly not exact on $E$,
hence $E$ is not IGCK. On the other hand,
a blow-up of a point in an LCK manifold
is again LCK (\cite[Proposition 2.4]{tric}, \cite[Theorem 1]{_Vuli_}), and the exceptional
divisor is IGCK, as well as all its submanifolds.

\hfill

\remark\label{_blow_up_Remark_} 
    {By contrast to the K\"ahler case, the
      class of manifolds of LCK type is not closed under
      blow-ups (\cite[Theorem 1.3, Claim 1.5]{ovv1}). Still, one obviously has that the blow-up
      of a manifold carrying a WLCK structure also carries
      a WLCK structure.}

\hfill

%
%

The next result, proven by Vaisman, states the dichotomy between K\"ahler and LCK manifolds:

\hfill

\theorem\label{_Vaisman_theorem_} (\cite[Theorem 2.1]{va_tr}, \cite[Section 4.3 for a different proof]{_OV_book_})\\ 
    {Let $(M,I,\omega,\theta)$ be a compact LCK manifold, $\dim_\C M\geq 2$. Assume that $(M, I)$ admits also a K\"ahler structure: then $[\theta]=0$, that is, $(M,I,g,\omega)$  is GCK.}
\endproof

\hfill

\remark\label{_WLCK_Vaisman_Remark_}(\cite[Lemma 2.5]{APV})
{If $(M,I,\omega,\theta)$ is assumed to be only WLCK, but still compact and carrying a K\"ahler metric, it follows similarly that $M$ must be GWLCK.}

\hfill

In our proofs, we shall use the following theorem: 

\hfill

\theorem\label{_fibrations_lemma_}(\cite[Lemma 3.1]{ovv1}, \cite[Lemma]{opv}) \\
Let $M$ be an LCK manifold, $B$ a  connected differentiable manifold,   $\dim B<\dim_\R M$, and
$\pi: M \rightarrow B$ a continuous, proper map. Assume that either:
\begin{description}
	\item[(i)] $B$ is an irreducible complex variety, and
	$\pi$ is holomorphic, or 
	\item[(ii)] $\pi$ is a locally trivial fibration with fibers complex subvarieties of $M$.
\end{description}
Suppose that 
any Lee class on $M$ is in the image of $\pi^*$ and the fibers of $\pi$ are positive  dimensional. Then any LCK structure on $M$ is GCK. \endproof

\hfill

\remark\label{_WLCK_fib_Lemma_Remark_}(\cite[Lemma 2.4]{APV})
{If $(M,I,\omega,\theta)$ is assumed to be only WLCK and such that the ``bad locus'' $\mathcal{B}(\omega)$ contains no fiber of $\pi$, it follows similarly that in fact $M$ is GWLCK.}

\hfill

In the next 3 subsections we provide more examples of LCK
manifolds. All known LCK manifolds fall whithin one of the
following three classes that we now describe.

\subsection{LCK manifolds with potential}\label{_LCK_pot_Subsection_}

\definition\label{_LCK_pot_} (\cite{ov_lckpot}, \cite[Chapter 12]{_OV_book_})\\ An LCK manifold has  {\bf LCK potential} if it
admits a K\"ahler covering on which the K\"ahler form
$\tilde \omega$
has a global and positive  potential function $\psi$,
$\tilde \omega= dd^c \psi$, such that
the deck group multiplies $\psi$ by a constant. 
In this case, $M$ is called {\bf an LCK manifold with
	potential}.

\hfill

\example Let $(M,I,\omega,g,\theta)$ be an LCK manifold with
the Lee form parallel with respect to the Levi-Civita
connection of $g$. Then $(M,I,\omega,g,\theta)$ is called {\bf a Vaisman manifold}. Let $\pi:(\tilde M,\tilde\omega)\rightarrow (M,\omega,\theta)$ be a  K\"ahler cover of a Vaisman manifold. Then one can see that the squared norm of $\pi^*\theta$ with respect to a K\"ahler metric $\tilde\omega$ is a global K\"ahler potential satisfying the conditions in \ref{_LCK_pot_}. Hence Vaisman manifolds are particular examples of LCK manifolds with potential. Among the examples of Vaisman manifolds we mention:
\begin{description}
	\item[(i)] All elliptic surfaces (\cite[Theorem 1]{bel}).
	\item[(ii)] The diagonal Hopf manifolds 
	$\frac{\C^n\backslash 0}{\langle A\rangle}$, where $A\in\GL(n,\C)$ is diagonalizable and its eigenvalues $a_i\in \C$ satisfy $0<|a_i|<1$ (\cite[Theorem 1]{go}, \cite[Section 2.5]{ov_pams}, \cite[Chapter 15]{_OV_book_}). 
\end{description}

\hfill

\remark\label{_canonical_foliation_remark_} Vaisman manifolds are endowed with a {\bf canonical foliation} locally generated by the Lee and anti-Lee fields $\theta^\sharp$ and $I\theta^\sharp$ (e. g. \cite[Theorem 3.1]{va_gd}). On a compact Vaisman manifold, any complex subvariety is tangent to the canonical foliation (\cite[Theorem 3.2]{ts}, \cite[Theorem 7.34]{_OV_book_}).

\hfill

The following theorem provides a useful criterion for a compact
LCK manifold to be of Vaisman type.

\hfill

\theorem (\cite[Proposition 3]{nico2})\label{_Killing_holo_Nicolina_Theorem_} 
Let $(M,I,\omega,\theta)$ be a compact LCK manifold, not GCK.
Consider a compact torus $T$ acting on $(M,I)$ by
biholomorphic diffeomorphisms.  Let $\goth t\subset TM$
be the Lie algebra of vector fields tangent
to this action. Assume that $I(\goth t) \cap \goth t\neq 0$.
Then $(M,I)$ is of Vaisman type. Moreover,
the Lie algebra $I(\goth t) \cap \goth t$
coincides with the Lie algebra generated by the
Lee and the anti-Lee fields: $\theta^\sharp, I\theta^\sharp$.  \endproof

\hfill

\remark  LCK manifolds with potential are stable to small deformations (\cite[Theorem 2.6]{ov_lckpot}). It follows that all {\bf linear Hopf manifolds} $\frac{\C^n\backslash 0}{\langle A\rangle}$, where $A\in\GL(n,\C)$, with eigenvalues $a_i\in \C$, $0<|a_i|<1$, are LCK with potential. On the other hand, linear but non-diagonal Hopf manifolds are not Vaisman (\cite[Theorem 2.16, Example 2.18]{ov_pams}).

\hfill

\theorem\label{_LCK_pot_embedded_Theorem_} 
(\cite[Theorem 3.4]{ov_lckpot}, \cite[Chapter 13]{_OV_book_})\label{_LCK_pot_embeds_in_Hopf_Theorem_}\\
 A compact LCK manifold with potential admits a
 holomorphic embedding into a linear Hopf
 manifold. Moreover, all compact LCK manifolds with
 potential contain an smooth elliptic curve. \endproof

\hfill

\remark Clearly,
any complex submanifold of an LCK manifold with potential
is LCK with potential. It was recently proven
(\cite[Corollary 4.3]{ov_non_linear}) that the non-linear Hopf manifolds
$\frac{\C^n\backslash 0}{\langle \gamma\rangle}$, where
$\gamma$ is an invertible, holomorphic contraction with
origin in $0\in\C^n$, can be holomorphically embedded into
linear Hopf manifolds; in particular, also the non-linear
Hopf manifolds are LCK with potential. 

\hfill

\remark\label{_cone_lee_classes_Remark_}
 One can easily prove that if some cohomology class
$[\theta]\in H^1(M,\R)$ is the class of the Lee form of some
LCK structure with potential $\psi$,  
then $u[\theta]$ is also the Lee class of an
LCK metric $\frac{dd^c \psi^u}{\psi^u}$
for any $u\in \R^{\geq 1}$. In fact, the set of Lee
classes on a compact LCK manifold with potential is an
open half-space in $H^1(M,\RR)$ (\cite[Theorem 8.4]{_OV_Lee_}).

\subsection{Manifolds of Inoue type}\label{_LCK_Inoue_type_Subsection_}

Recall that a real $(p,p)$-form $A$ on a complex
$n$-manifold $M$ is called {\bf weakly positive}
if $A \wedge \alpha^{n-p}$ is a non-negative
top form for any Hermitian form $\alpha$ on $M$.

\hfill

\definition\label{_defin_}
Let $A$ be a weakly positive, non-zero $(p,p)$-form
on a complex manifold $M$ admitting an LCK structure,
$\dim_\C M > p>0$. We say that $A$ {\bf consumes
the LCK structures} if for any LCK structure
$(\omega, \theta)$ on $M$, $\theta$ is cohomologous
to a closed 1-form $\theta_1$ such that $A \wedge \theta_1 =0$.
We will say that a compact complex manifold is {\bf of Inoue type} if it admits such an $A$ which is also closed, $dA=0$.

\hfill

\remark
As $A\not=0$ it follows that
the inequality $\int_M A\wedge \alpha^{n-p}\geq 0$ is strict,
\begin{equation}\label{int}
\int_M A\wedge \alpha^{n-p}>0
\end{equation}
for any Hermitian form $\alpha$.

\hfill

\example\label{_LCKOT_} 
First, we recall the LCK OT manifolds (see \cite{ot} for
details). Let $n\in \N_{>1}$ and fix a number field $K$
having $n-1$ real embeddindgs $\sigma_1,\dots,
\sigma_{n-1}$ and a single (up to conjugation) complex
one, say $\sigma_n.$ Let $\O_K$ be the ring of integers of
$K$ and $U\subset \O_{K}^*$ be a subgroup of finite index
formed by {\em positive units}, that is for any $u\in U$
one has $\sigma_1(u), ..., \sigma_{n-1}(u)>0.$ Consider the semidirect product
$\Gamma:=\O_K\rtimes U$ acting on
$$\H^{n-1}\times \C=\{(w_1,\dots, w_{n-1}, z)\ \vert \  w_i\in \H, z\in \C\}$$ 
as
$$a\cdot (w_1,\dots, w_{n-1}, z):=(w_1+\sigma_1(a),\dots, w_{n-1}+\sigma_{n-1}(a), z+\sigma_n(a))$$
for any $a\in \O_K$ and respectively
$$u\cdot (w_1,\dots, w_{n-1}, z):=(\sigma_1(u)w_1,\dots, \sigma_{n-1}(u)w_{n-1}, \sigma_n(u)z)$$
for any $u\in U.$ The  resulting manifold $M:=(\H^{n-1}\times \C)/\Gamma$ is a compact manifold that admits an LCK metric whose Lee form
is 
\begin{equation}\label{_leein_}
\theta_0:=d\log\left(\prod_{i=1}^{n-1} \Im(w_i)\right).
\end{equation}
Moreover, one can prove (\cite[Theorem 3.11]{_Otiman_}, \cite[Proposition 6.5]{_Istrati_Otiman_}) that for {\em any } LCK metric on $M$  the associated Lee form is cohomologous to the form $\theta_0$ defined above.

\hfill

\remark\label{_s0_}
Notice that if  in the previous example one takes $n=2$ one retrieves the LCK structure of the {\bf Inoue surfaces of type $S^0$} (for the description of these surfaces, see \cite{inoue}, also \cite{tric} and \cite[Chapters 22, 44]{_OV_book_}).

\hfill

\example\label{s+_}
Next we briefly recall the definition of Inoue surfaces of types $S^+$ and $S^{-}$ (\cite{inoue}, see also \cite[Chapter 44]{_OV_book_}). 
Namely, an {\bf Inoue surface of type $S^+$} is a quotient $M:=\H\times \C/\Gamma$ where $\Gamma$ is  a group of affine transformations generated by $g_0, g_1, g_2, g_3$ which are of the form
$$g_0(w, z):=(\alpha w, z+t)$$
and
$$g_i(w, z):=(w+a_i, z+b_iw+c_i), i=1,2,3$$
where 
$\alpha>1$ is a algebraic quadratic  integer, $t\in \C$ is some complex number and $a_i, b_i, c_i$ are appropiately chosen complex numbers.
It is known that $M$ has an LCK metric {\em if and only if} the parameter $t$ is {\em real} (\cite[Proposition 3.4]{tric}). Similarly to the previous case, the  associated Lee form of any LCK metric on $M$ is cohomologous to the Lee form given by \eqref{_leein_} (the result is implicit in \cite[Proof of Proposition 18]{bel}; see also \cite[Proposition 5.2]{_Apostolov_Dloussky_}).

The {\bf Inoue surfaces of type $S^{-}$} are quotients of
order $2$ of surfaces of type $S^+$ above {\em with real
  parameter $t$;} in particular, all of them  admit LCK
metrics and the de Rham class of the Lee form of any LCK metric on them is still unique, as in the previous cases.

\hfill

\proposition\label{_in are in_}
All LCK OT manifolds and all LCK Inoue  surfaces are of Inoue type, in the sense of \ref{_defin_}.

\hfill

\proof 
Let $M$ be any of the manifold as in the statement, and   $\H^{n-1}\times \C=\{(w_1,\dots, w_{n-1}, z)\ \vert \ w_i\in \H, z\in \C\}$
its universal cover. Define
$$A:=\left(\prod_{i=1}^{n-1}\frac{1}{y_i^2}\right)dx_1\wedge
dy_1\wedge\dots \wedge dx_{n-1}\wedge dy_{n-1}$$
where $w_i=x_i+\1 y_i.$
Clearly, $A$ is weakly positive, and satisfies $A \wedge \theta=0$,
where $\theta = d  \sum_i \log y_i$ is the standard LCK form.

On the other hand, keeping in mind the uniqueness of the Lee class of $\theta_0$ in \eqref{_leein_} we immediately see that $A$ consumes the LCK structures.
\endproof

\subsection{Manifolds with rational curves}
\label{_LCK with_rational_curves_Subsection_}
 
The ``Global Spherical Shell conjecture" (GSS conjecture)
claims that any minimal class VII surface $M$ with $b_2>0$
contains an open complex subvariety $U\subset M$
biholomorphic to a neighbourhood of the standard sphere
$S^3\subset \C^2$, and $M\backslash U$ is
connected. Surfaces which satisfy the GSS conjecture are
called {\bf Kato surfaces}. This conjecture is widely
believed to be true. Once the GSS conjecture is proven,
this finishes the classification of the compact complex
surfaces. If it is true, all non-K\"ahler surfaces are
LCK, except a particular class of $S^+$ Inoue
surfaces. See also \cite{_ovv:surf_} and \cite[Chapters
  24, 25]{_OV_book_} for an up to date treatment of LCK
geometry on compact complex surfaces.

The interesting feature about the remaining types of  LCK manifolds {\em known so far} is that they all carry rational curves.

An easy example of LCK manifold which is not with potential  is the blow-up of an LCK manifold with potential (obviously, the existence of a rational curve on the universal cover of the blow-up prevents the existence of a global plurisubharmonic function).
 
Recall that the blow-up $\hat{M}$ of an LCK manifold $M$
at a point (or, more generally, along an induced IGCK
submanifold)  is still LCK (\cite[Proposition 2.4]{tric}, \cite[Theorem 1]{_Vuli_}, \cite[Theorem 1.3]{ovv1}). Moreover,
it can be seen that for any Lee class $[\theta]\in
H^1(M)$, its pullback to $\hat{M}$ is the class of a Lee
form. 

\hfill

\example
A notable class of examples are the Kato manifolds (higher dimensional analogues of the Kato surfaces, see 
\cite{_Brunella:Kato_,iop,_IOPR:toric_Kato_}). Since they
also contain rational curves, it follows that they cannot
be LCK with potential. Indeed, let $M$ be an LCK manifold with potential and $C\subset M$ a rational curve. By \ref{_LCK_pot_embedded_Theorem_}, $M$ can be holomorphically embedded in a linear Hopf Manifold $H$, hence $C\subset H$, a contradiction because all curves on Hopf manifolds are elliptic.  

\hfill

\remark\label{_blow_up_of_Inoue_type_Remark_} One may see that the blow up at a point of an Inoue type LCK manifold is also of Inoue type and obviously contains a rational curve.

\section{Products with one factor LCK with potential are not LCK}


\lemma\label{_Vaisman_not_biholomo_to_product_Lemma_}
Let $M$ be a compact Vaisman manifold.
Then $M$ is not biholomorphic to a product of complex manifolds.

\hfill

\proof
Let $M= X \times Y$.
Since all positive-dimensional
subvarieties of a Vaisman manifold are tangent
to the canonical foliation (\ref{_canonical_foliation_remark_}), no positive-dimensional subvarieties
of a Vaisman manifold can intersect transversally.
Applying this to $X \times \{y\}$ ($y\in Y$),
and $\{x\}\times Y$ ($x\in X$), we arrive at contradiction.
\endproof

\hfill

This immediately brings:

\hfill

\proposition\label{_E_times_X_LCK_Corollary_}
Let $M=X \times E$ be the product of a compact
complex manifold and an elliptic curve. Then
$M$ does not admit a strict LCK structure.

\hfill

\proof
Consider the group of automorphisms of $E$
acting on $M$. Its tangent space is a Lie algebra ${\goth e}$
which satisfies ${\goth e}=I(\goth e)$. 
According to \ref{_Killing_holo_Nicolina_Theorem_},
$M$ is of Vaisman type. Now, by \ref{_Vaisman_not_biholomo_to_product_Lemma_},
we obtain that this is impossible.
\endproof

\hfill

The above result can be used to prove that any product of compact complex manifolds, one of which is LCK with potential, do not carry LCK metrics.

\hfill

\theorem\label{_nu_pot_}
Let $X$ be a compact LCK manifold with potential
and $Y$ any compact complex manifold, $\dim Y >0$.
Then the product $M:=X\times Y$ does not
admit a strict LCK structure.

\hfill

\proof
By \ref{_LCK_pot_embeds_in_Hopf_Theorem_}, $X$ contains an elliptic curve $E$.
A submanifold in an LCK manifold is also LCK.
Then, by \ref{_E_times_X_LCK_Corollary_}, the manifold
$E\times Y$ is globally conformally K\"ahler. Hence $\{x\}\times Y$ has an induced globally conformally K\"ahler (IGCK) structure for all $x\in X$.
We obtained that $M$ is fibered over $X$ with
IGCK fibers. This is impossible by \ref{_fibrations_lemma_}.
\endproof

\hfill

\corollary\label{_Nico_superseeded_Remark_}
The product of a compact
complex manifold $M$ with a compact complex curve admits no strict
LCK structure.

\hfill

\proof Recall that \cite[Theorem 7.8]{nico2} states that if a product of a compact LCK manifold $M$ with a compact complex curve has an LCK metric, then $M$ should be LCK with potential. By the above \ref{_nu_pot_}, such a situation cannot occur.

\section{Products with one factor of Inoue type are not LCK}

\theorem\label{_nu_inoue_}
Let $M= X \times Y$ be a product of two
manifolds admitting a strict LCK structure, and
$A$ a {closed} $(p, p)$-form on $X$ consuming the LCK structures.
Then $M$ does not admit an LCK structure.

\hfill

\pstep
By absurd, assume that $M$ admits an LCK structure
$(\omega, \theta)$. Replacing $\theta$ by a cohomologous
1-form $\theta'$, we can always change $\omega$ in its conformal
class resulting in an LCK structure $(\omega', \theta')$.

\hfill

{\bf Step 2:}
Let $\pi_1:\; M \arrow X$, $\pi_2:\; M \arrow Y$ be the projections.
Using the K\"unneth decomposition $H^1(M)= H^1(X) \oplus H^1(Y)$,
we may assume that $\theta$ is cohomologous
to $\pi_1^* \theta_1+ \pi_2^* \theta_2$, where $\theta_1, \theta_2$ are closed 1-forms on $X$, respectively $Y$.
Since the restriction of $(\omega, \theta)$
to $X\simeq X\times\{y\}= \pi_2^{-1}(y)$ is an LCK structure (for any $y\in Y$),
the cohomology class $[\theta_1]\in H^1(X)$
is a Lee class of a certain LCK structure.
Then $[\theta_1]$ contains 
a 1-form $\theta_1$ which satisfies $\theta_1\wedge A=0$.
Using Step 1, we can assume that
$\theta= \pi^*_1\theta_1 + \pi_2^* \theta_2$,
where $\theta_1\wedge A=0$.

\hfill

{\bf Step 3:} Let $n=\dim_\C X$.
Consider the form $B:= \pi_1^* A \wedge \omega^{n-p}$.
Recall that $\pi_1^*A \wedge\pi_1^* \theta_1=0$ and  $dA=0$. Then:
\begin{multline*}
d B= (n-p) \pi_1^*A \wedge \omega^{n-p} \wedge
(\pi_1^* \theta_1+ \pi_2^* \theta_2) \\=
(n-p) \pi_1^*A \wedge \omega^{n-p} \wedge\pi_2^* \theta_2
= (n-p) 
B \wedge\pi_2^* \theta_2 .
\end{multline*}
Consider the pushforward (that is, the fiberwise integral)
$(\pi_2)_* B\in C^\infty (Y)$ of $B$ to $Y$. Since $B$ is of type $(n, n)$ it follows that
$$(\pi_2)_* B(y)=\int_{\pi_2^{-1}(y)} \left(A\wedge \omega^{n-p}\right)\restrict{\pi_2^{-1}(y)}=\int_{\pi_2^{-1}(y)} A\wedge \left(\omega^{n-p}\right)\restrict{\pi_2^{-1}(y)}.$$  
Now, from the weak positivity of the $(p,p)$-form $A$, together with equation  \eqref{int}, we infer that the function $(\pi_2)_* B$ is positive and
nowhere vanishing. For all $\eta \in \Lambda ^* M$,  we have
\[
(\pi_2)_*(\eta \wedge \pi_2^* \theta_2) = (\pi_2)_*(\eta) \wedge \theta_2.
\]
Therefore
\[ d (\pi_2)_* B= (\pi_2)_* (dB) =  (n-p)(\pi_2)_* (B \wedge\pi_2^* \theta_2)
= (n-p)((\pi_2)_* B) \cdot \theta_2.
\]
This implies that $\theta_2 = \frac{1}{n-p}d\log ((\pi_2)_* B)$
is exact, contradicting the assumption that $Y$ is strictly LCK.
\endproof

\section{Products with one factor having  {rational curves} are not LCK}

\theorem\label{_nu_spl_}
Let $X$ be a compact strict LCK manifold and  $C\subset X$ a rational curve (that is, a closed analytic subspace whose normalisation is the projective line ${\mathbb P}^1$) 
and $Y$  a compact complex manifold, $\dim Y >0$.
Then the product $M:=X\times Y$ does not
admit a strict LCK structure.

\hfill

\proof
Let $(\omega, \theta)$ be an LCK structure on $M$ and $S$  the singular locus (possibly empty) of $C$. Blowing-up (possibly iterated) $M$  along $S\times Y$ we get a new manifold $\widehat{M}$ (which is in fact isomorphic to $\widehat{X}\times Y$, where $\widehat{X}$ is an embedded resolution of $C\subset X$). Let $E\subset \widehat{M}$ be the exceptional divisor of the blow-up $\sigma:\widehat{M}\ra M$ and let $\widehat{C}\subset \widehat{X}$ be the embedded resolution of $C\subset X$.
Then $\widehat{M}$ is WLCK, with structure given by $\left(\sigma^*(\omega), \sigma^*(\theta)\right)$; notice that ${\mathcal B}(\sigma^*(\omega))= E$.

Let $N:=\widehat{C}\times Y\subset \widehat{X}\times Y$; then $N$ is WLCK with the structure
$\left(\sigma^*(\omega)\restrict{N}, 
\sigma^*(\theta)\restrict{N}\right).$ Notice that the bad locus (see \ref{_WLCK_fib_Lemma_Remark_})   $\mathcal{B}(\sigma^*(\omega)\restrict{N})$ is just $N\cap E$.

The projection $pr_Y:N\ra Y$ has simply connected fibers, and thus induces an isomorphism  $H^1(Y)\ra H^1(N)$. It follows   that $\sigma^*(\theta)\restrict{N}$ is cohomologically a pull-back. Since $\mathcal{B}(\sigma^*(\omega)\restrict{N}$ is just $N\cap E$, hence intersecting the fibers of $pr_Y$ in finitely many points, we see that  \ref{_fibrations_lemma_} applies, hence $\sigma^*(\theta)\restrict{N}$ is cohomologically zero, henceforth $\theta\restrict{Y}$ is cohomologically trivial. Applying once again  \ref{_fibrations_lemma_} to $pr_X: X\times Y\ra X$ we get that $\theta$ is cohomologically trivial: this means that $M$ is GCK, hence  $X$ is GCK too, a contradiction, since we assumed $X$ to be strict LCK. \endproof


\hfill

\remark Owing to \ref{_Vaisman_theorem_}, the above \ref{_nu_spl_} and its proof also apply  to the case when instead of $C$ we consider any analytic subspace  $Z\subset X$ of  $\dim(Z)>1$ and such that  $Z$ has a desingularisation of K\"ahler type.

\section{Products of compact complex surfaces are not LCK}\label{_Proof_of_main_Section_}

\theorem\label{_No_LCK_on_products_of_surfaces_Theorem_}
Let $S$ be any compact complex surface. Assuming the GSS conjecture, then for any compact complex manifold $Y$, the product $M:=S\times Y$ has no LCK metric.

\hfill

\proof
By absurd, $M$ is LCK. Then, by \ref{_nu_spl_} $S$ must be minimal and not of Kato type.
Next, if  $S$ is elliptic or Hopf,  it would admit an LCK metric with potential, so this case is ruled out by \ref{_nu_pot_}.
Eventually, assuming true the GSS conjecture, we are left with the case when $S$ is an  Inoue surface, but this case is ruled out by \ref{_nu_inoue_}.
\endproof

\hfill

\noindent{\bf Acknowledgment:} We thank the anonymous referee for her or his extremely useful remarks.

{\scriptsize

{\small
	
	\noindent {\sc Liviu Ornea\\
		University of Bucharest, Faculty of Mathematics and Informatics, \\14
		Academiei str., 70109 Bucharest, Romania}, and:\\
	{\sc Institute of Mathematics ``Simion Stoilow" of the Romanian
		Academy,\\
		21, Calea Grivitei Str.
		010702-Bucharest, Romania\\
		\tt lornea@fmi.unibuc.ro,   liviu.ornea@imar.ro}
	
	\hfill

	\noindent {\sc Misha Verbitsky\\
		{\sc Instituto Nacional de Matem\'atica Pura e
			Aplicada (IMPA) \\ Estrada Dona Castorina, 110\\
			Jardim Bot\^anico, CEP 22460-320\\
			Rio de Janeiro, RJ - Brasil }\\
		also:\\
		Laboratory of Algebraic Geometry, \\
		Faculty of Mathematics, National Research University 
		HSE,\\
		6 Usacheva Str. Moscow, Russia}\\
	\tt  verbit@impa.br 

\hfill

	\noindent {\sc Victor Vuletescu\\
	University of Bucharest, Faculty of Mathematics and Informatics, \\14
	Academiei str., 70109 Bucharest, Romania}\\
	\tt vuli@fmi.unibuc.ro}}


\begin{thebibliography}{100}

\bibitem[APV]{APV} D. Angella, M. Parton, V. Vuletescu,
  {\em  On locally conformally K\"ahler threefolds with
    algebraic dimension two}, IMRN {\bf 5} (2023), 3948-3969. 

\bibitem[AD]{_Apostolov_Dloussky_} V. Apostolov, G. Dloussky, {\em On the Lee classes of locally conformally symplectic complex surfaces}, J. Sympl. Geom.  {\bf 16} (2018), 931-958.
	

\bibitem[Be]{bel} 
F. A. Belgun, {\em On the metric structure of non-K\"ahler complex surfaces}, Math. Ann. {\bf 317} (2000), 1-40.	




\bibitem[Br]{_Brunella:Kato_} M. Brunella, {\em Locally conformally K\"ahler metrics on Kato surfaces}, Nagoya Math. J. {\bf 202} (2011), 77-81.



\bibitem[GO]{go} P. Gauduchon, L. Ornea, {\em Locally conformally K\"ahler metrics on Hopf surfaces}, Ann. Inst. Fourier {\bf 48} (1998), 1107-1128.


\bibitem[In]{inoue} Ma. Inoue, {\em On surfaces of class $VII_0$}, Invent. Math., 24 (1974), 269-310.

\bibitem[Is]{nico2} 
N. Istrati, {\em Existence criteria
  for special locally conformally K\"ahler  metrics},
  Ann. Mat. Pura. Appl. {\bf 198} (2019), 335-353.
  
\bibitem[IO]{_Istrati_Otiman_}  N. Istrati, A. Otiman, {\em De Rham and twisted cohomology of Oeljeklaus-Toma manifolds}, Ann. Inst. Fourier (Grenoble) {\bf 69} (2019), no. 5, 2037-2066.

\bibitem[IOP]{iop}
N. Istrati, A.  Otiman, M. 
Pontecorvo, {\em On a class of  Kato manifolds},
IMRN, {\bf 7} (2021), 5366-5412. arXiv:1905.03224.

\bibitem[IOPR]{_IOPR:toric_Kato_}
N. Istrati, A. Otiman, M. 
 Pontecorvo, M.  Ruggiero {\em Toric Kato manifolds}, to appear in J. Ecole Polytechnique. arXiv:2010.14854 

%



\bibitem[OT]{ot} K. Oeljeklaus, M.  Toma, {\em Non-K\"ahler compact complex manifolds associated to number  fields}, Ann. Inst. Fourier {\bf 55}, no. 1 (2005), 1291-1300.

\bibitem[OPV]{opv} L. Ornea, M. Parton, V. Vuletescu, {\em Holomorphic submersions of locally conformally K\"ahler manifolds}, Ann. Mat. Pura Appl. {\bf (4) 193}  (2014), no. 5, 1345-1351.

\bibitem[OV1]{ov_lckpot} L. Ornea, M. Verbitsky, {\em Locally conformal K\"ahler manifolds with potential}, Math. Ann. {\bf 348} (2010), 25-33.

\bibitem[OV2]{ov_pams} L. Ornea, M. Verbitsky, {\em Locally conformally K\"ahler metrics obtained from pseudoconvex shells}, Proc. Amer. Math. Soc. {\bf 144} (2016), 325-335. 

\bibitem[OV3]{_OV_Lee_}  L. Ornea, M. Verbitsky, {\em Lee classes on LCK manifolds with potential}, arXiv:2112.03363. To appear in Tohoku Math. J.

\bibitem[OV4]{ov_non_linear} L. Ornea, M. Verbitsky, {\em Non linear Hopf manifolds are locally conformally K\"ahler}, J. Geom. Analysis {\bf 33} (2023) Article number: 201. arXiv:2202.12398.

\bibitem[OV5]{_OV_book_}  
L. Ornea, M. Verbitsky, Principles of locally conformally K\"ahler geometry, arXiv:2208.07188.

\bibitem[OVV]{ovv1}
 L. Ornea, M. Verbitsky, V. Vuletescu,
  {\em Blow-ups of locally conformally K\"ahler
    manifolds}, Int. Math. Res. Not. IMRN 2013, no. 12,
  2809-2821.





\bibitem[Ot]{_Otiman_} A. Otiman, {\em Morse-Novikov cohomology of locally conformally K\"ahler surfaces},   Math. Z. {\bf 289} (2018), no. 1-2, 605-628. arXiv:1609.07675.



\bibitem[Tr]{tric} 
F. Tricerri,  {\em Some examples of locally conformal K\"ahler manifolds}, 
Rend. Sem. Mat., Torino {\bf 40}, No.1  (1982), 81-92.

\bibitem[Ts1]{ts} K. Tsukada, {\em Holomorphic maps of compact
	generalized Hopf manifolds}, Geom. Dedicata {\bf 68} (1997), 61-71.

\bibitem[Ts2]{tsu} K. Tsukada, {\em The canonical foliation of a compact generalized Hopf manifold}, Differential Geom. Appl. {\bf 11} (1999), no. 1, 13-28. 


\bibitem[Va1]{va_tr} I. Vaisman, {\em On locally and globally conformal K\"ahler manifolds}, Trans. Amer. Math. Soc., {\bf 262} (1980), 533-542.

\bibitem[Va2]{va_gd} I. Vaisman, {\em Generalized Hopf manifolds}, Geom. Dedicata, {\bf 13} (1982), 231-255.


\bibitem[Vu]{_Vuli_} V. Vuletescu, {\em Blowing-up points on l.c.K. manifolds}, Bull. Math. Soc. Sci. Math. Roumanie (N.S.) {\bf 52(100)} (2009), no. 3, 387-390.



\bibitem[VVO]{_ovv:surf_} M. Verbitsky, V. Vuletescu, L. Ornea {\em Classification of non-K\"ahler surfaces and 
locally conformally K\"ahler geometry}, Russian Math. Surv. {\bf 76} (2021), 261-290. arxiv:1810.05768



\end{thebibliography}
\end{document}